\documentclass{amsart}

\swapnumbers \theoremstyle{plain}
\newtheorem{thm}{Theorem}[section]
\newtheorem{lem}[thm]{Lemma}

\newtheorem{prop}[thm]{Proposition}
\theoremstyle{definition}\newtheorem{rem}[thm]{Remark}
\newtheorem{defn}[thm]{Definition}

\newcommand{\Cal}{\mathcal}
\usepackage{amsfonts}
\newcommand{\D}{\mathbb{D}}
\newcommand{\R}{\mathbb{R}}

\newcommand{\C}{\mathbb{C}}
\newcommand{\N}{\mathbb{N}}

\newcommand{\g}{\overline\partial}
\DeclareMathOperator{\supp}{supp} 
 \DeclareMathOperator{\re}{Re}

 \numberwithin{equation}{section}

\begin{document}
\sloppy
\title{On the compactification of concave ends}
\author{Martin Brumberg and J\"urgen Leiterer}

\address{Martin Brumberg\\Institut f\"ur Mathematik \\Humboldt-Universit\"at zu Berlin
\\Rudower Chaussee 25\\D-12489 Berlin , Germany}
\email{martinbrumberg@googlemail.com}

\address{J\"urgen Leiterer\\Institut f\"ur Mathematik \\Humboldt-Universit\"at zu Berlin
\\Rudower Chaussee 25\\D-12489 Berlin , Germany}
\email{leiterer@mathematik.hu-berlin.de}

\thanks{2000 Mathematics Subject Classification. }
\maketitle

\section{Introduction}

\begin{defn}\label{3.8.09} We say
that $\rho:X\to\; ]a,b[$ is a {\bf $1$-corona} if:
\begin{itemize}
\item[]$a\in \R\cup\{-\infty\}$, $b\in \R\cup\{\infty\}$,
\item[]$X$ is a complex manifold,
\item[] $\rho$ is a strictly plurisubharmonic
$\Cal C^\infty$-function defined on $X$ and with values in the open interval $]a,b[$
such that the sets $\{\alpha\le \rho\le \beta\}$,
$a<\alpha<\beta<b$, are compacet.\end{itemize} We say that the
concave end of a $1$-corona $\rho: X\to]a,b[$ can be {\bf
compactified} if $X$ is (biholomorphic to) an open subset of a
complex space $\widehat X$ such that, for $a<c<b$, the set
$(\widehat X\setminus X)\cup \{a<\rho\le c\}$ is compact.
\end{defn}

 The concave end of a $1$-corona $\rho:X\to]a,b[$
always can be compactified if $n:=\dim X\ge 3$. This was proved by
Rossi \cite{Ro} and Andreotti-Siu \cite{AS}. For $n=2$ this not true
in general, as shown by a counterexample of Grauert, Andreotti-Siu
and Rossi \cite{AS,Gr,Ro}. However, if the concave end of a
$1$-corona $\rho:X\to]a,b[$ is even {\em hyperconcave} (i.e.
$a=-\infty$), then this is again true also for $\dim X=2$. This was
proved by Marinescu-Dinh \cite{MD}.

The Andreotti-Vesentini separation theorem \cite{AV} implies the
following necessary condition: If the concave end of a $1$-corona
$\rho:X\to]a,b[$ can be compactified, $\dim X\ge 2$, then $H^1(X)$
is Hausdorff, i.e. the space of exact $\Cal C^\infty_{0,1}$-forms is
closed with respect to uniform convergence on compact sets together
with all derivatives.

In the present paper we show that this condition is also sufficient.
We prove:
\begin{thm}\label{16.7.08}Let $\rho:X\to ]a,b[$ be a $1$-corona such that $H^1(X)$
is Hausdorff, and let $n:=\dim X\ge 2$.

Then $X$ can be embedded to $\C^{2n+1}$.

Hence, by the Harvey-Lawson theorem \cite{HaLa}, the  concave end of
$X$ can be compactified.
\end{thm}

\begin{rem} In view of the counterexample of Grauert,
Andreotti-Siu and Rossi mentioned above, theorem \ref{16.7.08} shows
that theorem 2 in \cite{Ra}, which states that $H^1(X)$ is Hausdorff
for each $1$-corona $\rho:X\to]a,b[$ also for $\dim X= 2$, cannot be
true.
\end{rem}

The present paper is a development of the PhD thesis of the first
author \cite{Br}, where the following theorem is proved:
\begin{thm}\label{5.8.08}
Let $\rho:X\to]a,b[$  ba a $1$-corona, $\dim X\ge 2$. Assume, for
some $a<t<b$, $H^1\big(\{a<\rho<t\}\big)$ is Hausdorff. Then there
exists  $\varepsilon>0$ such that  $\{t-\varepsilon<\rho<b\}$ can be
embedded to some $\C^N$.
\end{thm}

Theorem \ref{5.8.08} can be deduced as follows from theorem
\ref{16.7.08}: Suppose the hypotheses of theorem \ref{5.8.08} are
satisfied. Then, by theorem \ref{16.7.08},   $\{a<\rho<t\}$ can be
embedded to some $\C^N$. By the Harvey-Lawson theorem this implies
that the concave end of $\{a<\rho<t\}$ and, hence, the concave end
of $X$ can be compactified. By the Andreotti-Grauert theorem this
further implies that $H^1(X)$ is Hausdorff. Again using theorem
\ref{16.7.08}, we conclude that $X$ (and in particular
$\{t-\varepsilon<\rho<b\}$) can be embedded to some $\C^N$.  Note
however that in \cite{Br} a direct proof of theorem \ref{5.8.08} is
given, not using the compactification of the concave end.

The new contribution  of the present paper is the observation that
if
 $\rho:X\to]a,b[$ is a $1$-corona such that $H^1(X)$ is Hausdorff, then
 $H^1\big(\{a<\rho<t\}\big)$ is
Hausdorff for all $a<t<b$ (see theorem \ref{30.7.08'} below). After
that we could continue as in \cite{Br} to prove the following two
assertions on the existence of ''many'' holomorphíc functions:
\begin{itemize}
\item[(A)] The global holomorphic functions on $X$ separate points.
\item[(B)] For each $z\in X$ there exist $n$ global holomorphic
functions on $X$ which form a coordinate system at $z$.
\end{itemize}
From (A) and (B), the claim of theorem \ref{16.7.08} then follows in
the same way as in the proof of the embedding theorem for Stein
manifolds (see, e.g., theorem 5.3.6 in \cite{Ho}).

However, to prove (A) and (B) as in \cite{Br}, boundary estimates
for the $\g$-equation are needed, which are well-known but
nevertheless quite technical. To show that this is not necessary, in
the present paper we give another proof using only  inner estimates
for the $\g$-equation. A similar construction can be found already
in  \cite{HL2} (section 23), where assertions (A) and (B) are
deduced from certain stronger hypothesis, which implies the
Hausdorffness of $H^1\big(\{\rho<t\}\big)$ for each $a<t<b$.

{\bf Aknowledgement:}  We want to thank Genadi Henkin for the
discussions during the conference ''Complex Analyis and Geometry''
held from 23.06.08 to 29.06.08 in Bukarest (where the second autor
gave a lecture on the PhD-thesis of the first author). In
particular, we are grateful  for the hint that the methods proving
the embedability of Stein manifolds can be applied also to non-Stein
manifolds, provided properties (A) and (B) are already established.

\section{Each $H^1(\{\rho<t\})$ is Hausdorff if $H^1(X)$ is Hausdorff}

\begin{defn}\label{5.8.08''}
Let  $Y$ be a complex manifold of complex dimension $n$.

Then we denote by $\Cal E^{p,q}(Y)$, $0\le p,q\le n$, the Fr\'echet
space of $\Cal C^\infty_{p,q}$-forms on $Y$ endowed with the
topology of uniform convergence together with all derivatives on the
compact subsets of $Y$.

If $K\subseteq Y$ is a compact set, then we denote by $\Cal
D^{p,q}_K(Y)$ the subspace of $\Cal E^{p,q}(Y)$ which consists of
the forms with support in $K$, endowed with the same topology.

By $\Cal D^{p,q}(Y)$ we denote the space all $\Cal
C^\infty_{p,q}$-forms on $Y$ with compact support, endowed with the
Schwartz topology, i.e. the finest local convex topology such that,
for each compact set $K\subseteq Y$, the embedding $\Cal
D^{p,q}_K(Y)\to\Cal D^{p,q}(Y)$ is continuous.
\end{defn}

Recall  the following well known fact from Serre duality \cite{L}:

\begin{prop}\label{30.7.08} Let $Y$ be a complex manifold of dimension $n$. Then
the space $\g \Cal E^{0,0}(Y)$ is closed in $\Cal E^{0,1}(Y)$, if
and only if, $\g\Cal D^{n,n-1}(Y)$ is closed in $\Cal D^{n,n}(Y)$.
\end{prop}

We need also the following supplement to this (see theorem 2.7. in
\cite{LL}):

\begin{prop}\label{6.8.08} Let $Y$ be a complex manifold of dimension $n$.
Then the space $\g \Cal D^{n,n-1}(Y)$ is closed in $\Cal
D^{n,n}(Y)$, if and only if, for each compact set $K\subseteq Y$,
the space $ \Cal D^{n,n}_K(Y)\cap\g\Cal D^{n,n-1}(Y)$ is closed in
$\Cal D^{n,n}_K(Y)$.
\end{prop}

Moreover, we shall use the following result from Andreotti-Grauert
theory (see, e.g., theorem 16.1 in \cite{HL2}):

\begin{prop}\label{7.8.08} Let $\rho:X\to ]a,b[$ be a $1$-corona,
and let $a<c<b$. Then each $\g$-closed $f\in \Cal E^{n,n-1}(X)$ with
$f\equiv 0$ on $\{c<\rho<b\}$ is $\g$-exact.
\end{prop}

Using these propositions, we now obtain:

\begin{thm}\label{30.7.08'} Let $\rho:X\to]a,b[$ be a $1$-corona
such that $H^1(X)$ is Hausdorff. Then $H^1\big(\{\rho<t\}\big)$ is
Hausdorff for all $a<t<b$.
\end{thm}

\begin{proof} Let $a<t<b$, and let $K\subseteq \{\rho<t\}$ be
a compact set. By propositions \ref{30.7.08} and \ref{6.8.08}, it is
sufficient to prove that the space $ \Cal
D^{n,n}_K\big(\{\rho<t\}\big)\cap\g\Cal
D^{n,n-1}\big(\{\rho<t\}\big)$ is closed in the Fr\'echet space
$\Cal D^{n,n}_K\big(\{\rho<t\}\big)$.

Let a sequence $f_\nu\in \Cal
D^{n,n}_K\big(\{\rho<t\}\big)\cap\g\Cal
D^{n,n-1}\big(\{\rho<t\}\big)$ be given which converges to some
$f\in \Cal D^{n,n}_K\big(\{\rho<t\}\big)$ with respect to the
topology of $\Cal D^{n,n}_K\big(\{\rho<t\}\big)$. We have to find
$w\in \Cal D^{n,n-1}\big(\{\rho<t\}\big)$ with $\g w=f$.

Extending by zero, we may view $f_\nu$ as forms in $\g\Cal
D^{n,n-1}(X)$, and $f$ as a form in $\Cal D^{n,n}(X)$. As $H^1(X)$
is Hausdorff it follows from proposition \ref{30.7.08} that $\g\Cal
D^{n,n-1}(X)$ is closed in $\Cal D^{n,n}(X)$ with respect to the
Schwartz topology. Moreover, as all $f_\nu$ belong to $\g\Cal
D^{n,n-1}(X)$ and the embedding $\Cal D^{n,n}_K(X)\to \Cal
D^{n,n}(X)$ is continuous (by definition of the Schwartz topology),
$f_\nu$ converges to $f$ with respect to the Schwartz topology of
$\Cal D^{n,n}(X)$. Hence $f\in \g\Cal D^{n,n-1}(X)$, i.e. $f=\g u$
for some $u\in \Cal D^{n,n-1}(X)$.

Choose $a<s_f<s'_f<t$ with  $\supp f\subseteq K\subseteq
\{s_f<\rho<s'_f\}$. Moreover, as $u$ has compact support in $X$, we
can find $a<s_u<s_f<t<s_u'<b$ such that $\supp
u\subseteq\{s_u<\rho<s'_u\}$.  Then $\g u=f=0$ on $\{s_f'<\rho\}$,
and, since $\rho:\{s_f'<\rho\}\to ]s_f',b[$ is a $1$-corona, this
implies by proposition \ref{7.8.08}
 that $u\big\vert_{\{s_f'<\rho\}}=\g v$ for some $v\in \Cal
E^{n,n-2}\big(\{s_f'<\rho\}\big)$.  Choose $s'_f<t'<t$ and a $\Cal
C^\infty$-function $\chi$ on $X$ such that $\chi\equiv 1$ on
$\{t'\le\rho\}$ and $\chi\equiv 0$ on $\{\rho\le s_f'\}$. Set
\[w=\begin{cases}u\qquad&\text{on }\{\rho\le t'\},\\
u- \g (\chi v)\qquad&\text{on }\{t'<\rho\}.
\end{cases}\] Then $w\in \Cal D^{n,n-1}_{\{s'_u\le\rho\le t'\}}(X)$ and $\g w=f$ on $X$.
Since $\{s'_u\le\rho\le t'\}$ is a compact subset of $\{\rho<t\}$,
this completes the proof.
\end{proof}

\section{Approximation}

Here we prove the following theorem:

\begin{thm}\label{30.7.08-} Let $\rho:X\rightarrow ]a,b[$ be
a $1$-corona such that $H^1(X)$ is Hausdorff, and let $a<t<b$. Then each
holomorphic function defined in a neighborhood of $\{\rho\le t\}$
can be approximated uniformly on $\{\rho\le t\}$ by holomorphic
functions defined on $X$.
\end{thm}

We prove this by means of Grauert's bump method, similarly as in the
case of strictly pseudoconvex manifolds (see, e.g., the proof of
theorem 2.12.3 in \cite{HL1}), although, in distinction to strictly
pseudoconvex manifolds, here the $\g$-equation not always can be
solved. This is possible, because for the bump method only the
Hausdorffness (and not the vanishing) of all
$H^1\big(\{\rho<t\}\big)$, $a<t<b$, is needed.

\begin{defn}\label{22.7.08} A collection $$\Big(\rho:X\to]a,b[\;,\;
\varphi:X\to]a,b[\;,\;c\;,\;(U;z_1,\ldots,z_n)\Big)$$ will be called
a {\bf Grauert bump} if:

\begin{itemize}
\item[(a)] $\rho:X\to]a,b[$ and
$\varphi:X\to]a,b[$ are strictly pseudoconvex coronas\footnote{with
the same $X$ and $]a,b[$ but different $\rho$ and $\varphi$} such
that $\{\rho\le c\}\subseteq\{\varphi\le c\}$;
\item[(b)] $U$ is a relatively compact open subset of $X$, and $z_1,\ldots,z_n$
are holomorphic coordinates defined in a neighborhood of $\overline
U$ such that $U$ is a ball with respect to these coordinates and $
\rho=\varphi$ in a neighborhood of $X\setminus U$.
\end{itemize}
\end{defn}

\begin{lem}\label{22.7.08-} Let
$$\Big(\rho:X\to]a,b[\;,\;
\varphi:X\to]a,b[\;,\;c\;,\;(U;z_1,\ldots,z_n)\Big)$$ be a Grauert
bump. Then each holomorphic function defined in a neighborhood of
$\{\rho\le c\}$ can be approximated uniformly on $\{\rho\le c\}$ by
holomorphic functions defined in a neighborhood of $\{\varphi\le
c\}$.
\end{lem}

\begin{proof} Let a holomorphic function $f$ in a neighborhood of $\{\rho\le c\}$
be given. Choose $\varepsilon>0$ so small that
$$
\Big(\rho:X\to]a,b[\;,\;
\varphi:X\to]a,b[\;,\;c+\varepsilon,\;(U;z_1,\ldots,z_n)\Big)$$ is
still a Grauert bump and $f$ is holomorphic in a neighborhood of
$\{\rho\le c+\varepsilon\}$. As $\overline U$ is convex with respect
to the coordinates $z_1,\ldots,z_n$ and the function $\rho$ is
plurisubharmonic in a neighborhood of $\overline U$, the set $
\overline U\cap\{\rho\le c+\varepsilon\}$ is polynomially convex
with respect to $z_1,\ldots,z_n$ (see, e.g., theorem 4.3.2 in
\cite{Ho} or theorem 2.7.1. in \cite{HL1}). Therefore we can find a
sequence of holomorphic functions $u_\nu$ defined on a neighborhood
of $\overline U$, which converges to $f$ uniformly on $\overline
U\cap \{\rho\le c+\varepsilon\}$.

Since $\rho=\varphi$ in a neighborhood of $X\setminus U$, we can
find a relative compact open subset $U'$ of $U$ such that
$\rho=\varphi$ outside $U'$. Choose a $\Cal C^\infty$-function
$\chi$ on $X$ such that $\chi\equiv 1$ in $U'$ and  $\chi\equiv 0$
in a neighborhood of $X\setminus U$. Setting
\[
v_\nu=\begin{cases}f+\chi (u_\nu-f)\quad&\text{on }U\cap\{\rho\le c+\varepsilon\},\\
f&\text{on }(X\setminus U')\cap\{\varphi\le
c+\varepsilon\}=(X\setminus U')\cap\{\rho\le c+\varepsilon\},\\
u_\nu&\text{on }U',
\end{cases}\]we obtain a sequence of $\Cal C^\infty$-functions $v_\nu$ defined
in a neighborhood of
 $\{\varphi\le c+\varepsilon\}$ such that
\begin{equation}\label{3.8.08-}
v_\nu\text{ converges to }f\text{ uniformly on }\{\rho\le
c+\varepsilon\},
\end{equation}
\begin{equation}\label{3.8.08--}
\g v_\nu\text{ converges to zero in the Fr\'echet topology of }\Cal
E^{0,1}\big(\{\varphi<c+\varepsilon\}\big).
\end{equation}Choose $t$ with $a<t<c$ so close to $a$ that
$\overline U\subseteq\{\varphi>t\}=\{\rho>t\}$. Since, by theorem
\ref{30.7.08'}, $\g \Cal E^{0,0}(\{\varphi<c+\varepsilon\}$ is
closed with respect to the Fr\'echet topology of  $\Cal
E^{0,1}\big(\{\varphi<c+\varepsilon\}\big)$, then, by the Banach
open mapping theorem, it follows from \eqref{3.8.08--} that there
exists a sequence of functions $w_\nu\in \Cal
E^{0,0}\big(\{\varphi<c+\varepsilon\}\big)$ such that
\begin{align}\label{4.8.08}&\g w_\nu=\g v_\nu\text{ on }\{\varphi<c+\varepsilon\}\\
\label{4.8.08''}&\text{and } w_\nu \text{ converges to zero
uniformly on } \left\{t\le\varphi\le
c+\frac{\varepsilon}{2}\right\}.
\end{align}
Then, in a neighborhood of $\{\varphi\le t\}$,
$$
\g w_\nu=\g v_\nu=(u_\nu-f)\g \chi=0.
$$Therefore, the functions  $w_\nu$ are holomorphic in a neighborhood of $\{\varphi\le
t\}$.

Since, for each sufficiently small $\delta>0$,
$\varphi:\{\varphi<t+\delta\}\to]a,t+\delta[$ is a $1$-corona, this
implies that
\begin{equation}\label{4.8.08'}
\sup\limits_{ \varphi(\zeta)\le t}\big\vert
w_\nu(\zeta)\big\vert=\max\limits_{\varphi(\zeta)=t}\big\vert
w_\nu(\zeta)\big\vert.
\end{equation}Indeed,  by the lemma of Morse
(see, e.g., proposition 0.5 in Appendix B of \cite{HL2}), we may
assume that $\varphi$ has only non-degenerate critical points.  Now
first let $a<s'<s''\le t$ such that no critical point of $\varphi$
lies on $\{s'\le\varphi\le s''\}$. Then, for $s'\le s\le s''$,
 locally, the surface $\{\varphi= s\}$ is strictly convex with
respect to appropriate  holomorphic coordinates (as the boundary of
$\{\varphi<s\}$). Hence, for each point $\xi\in \{\varphi<s\}$
sufficiently close to $\{\varphi=s\}$, there is a smooth complex
curve $L\subseteq X$  such that $L\cap\{\varphi\le s\}$ is compact
and $\xi\in L\cap \{\varphi<s\}$, which yields
$$w_\nu(\xi)\le \max\limits_{\zeta\in L\cap\{\varphi=s\}}\big\vert
w_\nu(\zeta)\big\vert\le \max\limits_{\varphi(\zeta)=s}\big\vert
w_\nu(\zeta)\big\vert.$$ Therfore: whenever $a<s'<s''\le t$ such
that  $\varphi$ has no critical points on $\{s'\le\varphi\le s''\}$,
then
$$
\max\limits_{s'\le \varphi(\zeta)\le s''}\big\vert
w_\nu(\zeta)\big\vert=\max\limits_{\varphi(\zeta)=s''}\big\vert
w_\nu(\zeta)\big\vert.
$$
As the critical points of $\varphi$ are isolated, now
\eqref{4.8.08'} now follows by continuity.

It remains to set
$$
h_\nu=v_\nu-w_\nu.
$$
By \eqref{4.8.08} these functions $h_\nu$ are holomorphic on a
neighborhood of  $\{\varphi\le c\}$, and by \eqref{4.8.08''} and
\eqref{3.8.08-}, they converge to $f$ uniformly on  $\{\rho\le c\}$.

\end{proof}

Theorem \ref{30.7.08-} now is an immediate consequence of the
preceding lemma and the following lemma:

\begin{lem}\label{22.7.08'}
Let $\rho:X\to ]a,b[$ be a strictly pseudoconvex corona, and let
$a<t<b$. Then there exists a sequence
$$\Big(\rho_j:X\to]a,b[\;,\;
\varphi_j:X\to]a,b[\;,\;c_j\;,\;
(U_j;z_1^{(j)},\ldots,z_n^{(j)})\Big)\,,\qquad j\in \N,$$ of Grauert
bumps such that
\begin{itemize}
\item $\rho_0=\rho$ and $c_0=t$,
\item $\rho_{j+1}=\varphi_j$ for all $j\in \N$,
\item $X=\bigcup_{j=0}^\infty \{a<\rho_j\le c_j\}$.
\end{itemize}
\end{lem}

The latter lemma follows easily from the fact that a $\Cal
C^2$-small pertubation of a strictly plurisubharmonic function is
again such a function. We omit the proof, which is the same as in
the case of a strictly pseudoconvex manifold (see, e.g., lemma
2.12.4 in \cite{HL1} or lemma 12.3 in \cite{HL2} and the subsequent
remark).

\section{Proof of theorem \ref{16.7.08}}

In this section $\rho:X\to ]a,b[$ is a $1$-corona such that $H^1(X)$
is Hausdorff. As already observed in the introduction, we only have
to prove that conditions (A) and (B) are satisfied. To prove this,
let two different points $\xi,\eta\in X$ and a holomorphic
tangential vector $\theta$ of $X$ at $\xi$  be given. It is
sufficient to show that:
\begin{itemize}
\item[(A')] There exists a holomorphic function $f$ on $X$ such
that $f(\xi)\not=f(\eta)$.
\item[(B')] There exists a sequence of holomorphic functions $g_\nu$
on $X$ such that $\partial g_\nu(\xi)$ converges to $\theta$.
\end{itemize}

We may assume that $\rho(\xi)\ge \rho(\eta)$. Choose a system
$z=(z_1,\ldots,z_n)$ of holomorphic coordinates defined in a
neighborhood $U$ of $\xi$ such that $\eta\not\in U$, $z(\xi)=0$ and
$z(U)$ contains the unit ball. Let
$$ F:=-2\sum_{j=1}^n\frac{\partial\rho(\xi)}{\partial
z_j}\, z_j-\sum_{j,k=1}^n\frac{\partial^2\rho(\xi)}{\partial
z_j\partial z_k}\,z_jz_k
$$
be the Levi polynomial with respect to $z_1,\ldots,z_n$ of $\rho$ at
$\xi$. As well known (see, e.g., \cite{HL1}), then we can find
$0<\varepsilon<1/2$ with
$$
\re F\ge \rho(\xi)-\rho+2\varepsilon\vert z\vert^2\qquad\text{on
}\big\{\vert z\vert\le 2\varepsilon\big\}.
$$
It follows that
$$
\re F\ge \varepsilon^3\qquad\text{on
}\big\{\rho\le\rho(\xi)+\varepsilon^3\big\}\cap
\big\{\varepsilon\le\big\vert z\big\vert\le 2\varepsilon\big\}.
$$
Set $h=e^{-F}$ on $U$. Then
\begin{align*}
&h(\xi)=1\text{ and}\\
 &\vert h\vert\le e^{-\varepsilon^3}<1\text{ on }\big\{
\rho\le\rho(\xi)+\varepsilon^3\big\}\cap\big\{\varepsilon\le \vert
z\vert\le 2\varepsilon\big\}.
\end{align*}
Choose a $\Cal C^\infty$-function $\chi$ on $X]$ with $\chi\equiv 1$
in a neighborhood of  $\big\{\vert z\vert\le \varepsilon\}$ and
$\chi\equiv 0$ in a neighborhood of $X\setminus\big\{\vert
z\vert<2\varepsilon\big\}$. Setting
\[
f'_\nu=\begin{cases}h^\nu\chi\qquad&\text{on }\big\{\vert z\vert
<2\varepsilon\big\},\\
0&\text{on }X\setminus \big\{\vert
z\vert<2\varepsilon\big\},\end{cases}
\]we obtain a sequence $f'_\nu$ of $\Cal C^\infty$-functions on $X$ with
\begin{align}\label{8.8.08}
&f'_\nu(\xi)=1\text{ for all }\nu,\\
\label{8.8.08'}
&f'_\nu(\eta)=0\text{ for all }\nu\text{ and}\\
\label{8.8.08''}&\g f'_\nu=h^\nu\g\chi\text{ converges to zero
uniformly on }\big\{ \rho\le\rho(\xi)+\varepsilon^3\big\}.
\end{align}Moreover, let $\theta_j$ be the coefficients with
$\theta=\sum_{j=1}^n\theta_j\frac{\partial}{\partial
z_j}\big\vert_\xi$. Then, setting
\[ g'_\nu=\begin{cases}\Big(\sum_{j=1}^n\theta_jz_j\Big)
f'_\nu\qquad&\text{on }\{\vert z\vert<2\varepsilon\},\\0&\text{on
}X\setminus \{\vert z\vert<2\varepsilon\},\qquad
\end{cases}
\]we obtain a sequence $g'_\nu$ of $\Cal
C^\infty$-functions on $X$ with
\begin{align}\label{8.8.08-}
&\partial g'_\nu(\xi)=
f'_\nu(\xi)\,\theta+\bigg(\sum_{j=1}^n\theta_jz_j(\xi)\bigg)\partial
f'_\nu(\xi)=\theta\text{ for all }\nu\text{ and}\\
\label{8.8.08--}&\g g'_\nu=\bigg(\sum_{j=1}^n\theta_j z_j\bigg)\g
f'_\nu\text{ converges to zero uniformly on }\big\{
\rho\le\rho(\xi)+\varepsilon^3\big\}.
\end{align}
Since, by theorem \ref{30.7.08'}, $\g \Cal
E^{0,0}\big(\{\rho<\rho(\xi)+\varepsilon^3\}\big)$ is closed in the
Fr\'echet space $\Cal
E^{0,1}\big(\{\rho<\rho(\xi)+\varepsilon^3\}\big)$, now, by
\eqref{8.8.08''}, \eqref{8.8.08--} and the Banach open mapping
theorem, we can find  sequences $f''_\nu, g''_\nu\in \Cal
E^{0,0}\big(\{\rho<\rho(\xi)+\varepsilon^3\}\big)$ which converge to
zero uniformly on $\{\rho\le \rho(\xi)+\varepsilon^3/2\}$ such that
$$\g f''_\nu=\g f'_\nu\qquad{and}\qquad  \g g''_\nu=\g g'_\nu\qquad\text{for all }\nu.
$$
Setting $f'''_\nu=f'_\nu-f''_\nu$ and $g'''_\nu=g'_\nu-g''_\nu$ we
get holmorphic functions on $\{\rho<\rho(\xi)+\varepsilon^3\}$ such
that, by \eqref{8.8.08}, \eqref{8.8.08'} and \eqref{8.8.08-},
\begin{align}\label{8.8.08*}
&f'''_\nu(\xi)\text{ converges to }1,\\
\label{8.8.08**}&f'''_\nu(\eta)\text{ converges to }0\text{ and}\\
\label{8.8.08***} &\partial g'''_\nu(\xi)\text{ converges to }
\theta.
\end{align}

{\em Proof of (A'):} By \eqref{8.8.08*} and \eqref{8.8.08**} we can
find $\nu_0$ so large that $f'''_{\nu_0}(\xi)>3/4$ and
$f'''_{\nu_0}(\eta)<1/4$. By theorem \ref{30.7.08-}, now we can find
a holomorphic function $f$ on $X$ such that $\vert f-f'''\vert<1/4$
on $\{\rho\le \rho(\xi)\}$. As $\xi,\eta\in \{\rho\le \rho(\xi)\}$,
then $f(\xi)>1/2>f(\eta)$.

{\em Proof of (B'):} Again by theorem \ref{30.7.08-}, we can find a
sequence of holomorphic functions $g_\nu$ on $X$ such that $g_\nu-
g'''_\nu$ converges to zero uniformly on $\big\{\rho\le
\rho(\xi)+\varepsilon^3/2\big\}$. Since $\xi$ is an inner point of
$\big\{\rho\le \rho(\xi)+\varepsilon^3/2\big\}$, then also $\partial
g_\nu(\xi)-\partial g'''_\nu(\xi)$ converges to zero. By
\eqref{8.8.08***} this implies that $\partial g_\nu(\xi)$ converges
to $\theta$.


\begin{thebibliography}{GuRo}

\bibitem[AS]{AS} A. Andreotti, Y.-T. Siu, {\em Projective embeddings
of pseudoconcave spaces,} Ann. Sc. Norm. Super. Pisa {\bf 24},
231–278 (1970).

\bibitem[AV]{AV} A. Andreotti, E. Vesentini, {\em Carleman estimates for the
Laplace-Beltrami equation on complex manifolds,} Publ. Math., Inst.
Hautes Etud. Sci. {\bf 25}, 81-130 (1965).

\bibitem[Br]{Br} M. Brumberg, {\em Kompaktifizierung streng
pseudokonkaver Enden und Separiertheit der ersten
Dolbeaul-Kohomologie,} Dissertation am Institut f\"ur Mathematik der
Humboldt-Universit\"at zu Berlin, verteidigt am 8.5.2008.

\bibitem[Gr]{Gr} H. Grauert, {\em Theory of q-convexity and q-concavity,}
 Several Complex Variables VII, H. Grauert, Th. Peternell, R. Remmert, eds. Encyclopedia of
Mathematical Sciences, vol. {\bf 74}, Springer 1994.

\bibitem[HaLa]{HaLa} R. Harvey, B. Lawson, {\em On boundaries of complex
analytic varieties. I,} Ann. Math. (2) {\bf 102}, 223-290 (1975).

\bibitem[HL1]{HL1} G. Henkin, J. Leiterer, {\em Theory of functions on complex manifolds},
Birkh\"auser 1984.

\bibitem[HL2]{HL2} G. Henkin, J. Leiterer, {\em Andreotti-Grauert theory by integral formulas},
Progress in Mathematics {\bf 74}, Birkh\"auser 1988.

\bibitem[Ho]{Ho} L. H\"ormander, {\em An introduction
to complex analysis in several variables,} North-Holland 1990.

\bibitem[L]{L} H. Laufer, {\em On Serre duality and envelopes of
holomorphy}, Trans. Amer. Math. Soc. {\bf 128},414-436 (1967).

\bibitem[LL]{LL} C. Laurent-Thi\'ebaut, J. Leiterer, {On Serre duality},
Bull. Sci. Math. {\bf 124} , 93-106 (2000).

\bibitem[MD]{MD} G. Marinescu, T.-C. Dinh, {\em On the
compactification of hyperconcave ends and the theorems of Siu-Yau
and Nadel,} Invent. math. {\bf 164}, No. 2, 233-248 (2006).

\bibitem[Ra]{Ra} Ramis, J.-P., {\em Th\'eor\`emes de s\'eparation et de
finitude pour l'homologie et la cohomologie des espaces
$(p,q)$-convexes-concaves,} Ann. Scuola Norm. Sup. Pisa Cl. Sci {\bf
27}, 933-997 (1973).

\bibitem[Ro]{Ro} Rossi, H., {\em Attaching analytic spaces to an
analytic space along a pseudoconcave boundary}, Proc. Conf. Complex.
Manifolds (Minneapolis), pp. 242–256,  Springer 1965.




\end{thebibliography}
\end{document}